\documentstyle[12pt,amssymb,amstex]{amsart}
\numberwithin{equation}{section}

\def\bh{{\mathbb H}}

\def\bn{{\mathbb N}}

\def\br{{\mathbb R}}

\def\Q{\mathbf{Q}}
\def\H{\mathbf{H}}
\def\R{\mathbf{R}}
\def\P{\mathbf{P}}
\newtheorem{thm}{Theorem}[section]
\newtheorem{lem}[thm]{Lemma}
\newtheorem{cor}[thm]{Corollary}
\newtheorem{prop}[thm]{Proposition}

\begin{document}

\title[On the ergodic principle]
{On the ergodic principle for markov and quadratic Stochastic
Processes and its relations}
\author{Nasir Ganikhodjaev}
\address{Nasir Ganikhodjaev\\
Department of Mechanics and Mathematics\\
National University of Uzbekistan\\
Vuzgorodok, 700174, Tashkent, Uzbekistan\\
and Faculty of Science,\\
IIUM, 53100 Kuala Lumpur, Malaysia} \email{{\tt
nasirgani@@yandex.ru}}
\author{Hasan Akin}
\address{Hasan Akin\\
Department of Mathematics\\
Arts and Science Faculty\\
Harran University, Sanliurfa, 63200, Turkey} \email{{\tt
akinhasan@@harran.edu.tr}}
\author{Farrukh Mukhamedov}
\address{Farruh Mukhamedov\\
Department of Mechanics and Mathematics\\
National University of Uzbekistan\\
Vuzgorodok, 700174, Tashkent, Uzbekistan} \email{{\tt
far75m@@yandex.ru}}

\begin{abstract}
In the paper we prove that  a quadratic stochastic process
satisfies the ergodic principle if and only if the associated
Markov process satisfies one.

\vskip 0.3cm \noindent {\bf Mathematics Subject Classification}:
60K35, 60J05, 60F99, 92E99, 47A35.\\
{\bf Key words}: Quadratic stochastic processes; Markov processes,
Ergodic principle.
\end{abstract}

\maketitle

\section{Introduction}

It is known that Markov processes are well-developed field of
mathematics which have various applications in physics, biology
and so on. But there are some physical models which cannot be
described by such processes. One of such models is a model related
to population genetics. Namely,  this model is described by
quadratic stochastic processes (see \cite{L} for review). To
define it, we denote
\begin{eqnarray*}
& & \ell^{1}=  \{x=(x_{n}):\|
x\|_{1}=\overset{\infty}{\underset{n=1}{\sum}}|x_{n}|<\infty ;\
x_{n}\in \mathbb{R}\},\\
&& S =\{x\in \ell^{1}:x_{n}\geq  0 ;\ \|x\|_{1}=1\}.
\end{eqnarray*}
Hence this process is defined as follows (see
\cite{G1},\cite{SG1}): Consider a family of functions $\{
P_{ij,k}^{[s,t]}: i,j,k\in\bn,\ s,\ t\in \br^{+},\ t-s \geq 1\}$.
This family is said to be {\it quadratic stochastic process
(q.s.p.)} if for fixed $s,t\in\br_+$
it satisfies the following conditions:\\
\begin{enumerate}
   \item[(i)] $P_{ij,k}^{[s,t]}=P_{ji,k}^{[s,t]}$ for any $i,j,k\in\mathbb{N}$.
   \item[(ii)] $P_{ij,k}^{[s,t]}\geq 0$ and
$\overset{\infty}{\underset{k=1}{\sum }}P_{ij,k}^{[s,t]}=1$ for
any $i,j,k\in \mathbb{N}$.
   \item[(iii)] An analogue of Kolmogorov-Chapman equation; here there are two
variants: for the initial point $x^{(0)}\in S$,
$x^{(0)}=(x^{(0)}_{1},x^{(0)}_{2},\cdots)$ and $s<r<t$ such that $t-r\geq 1$,$r-s\geq 1$\\
(iii$_{A}$)
$$
P_{ij,k}^{[s,t]}=\overset{\infty}{\underset{m,l=1}{\sum
}}P_{ij,m}^{[s,r]}P_{ml,k}^{[r,t]}x^{(r)}_{l}, $$ where
$x^{(r)}_{k}$ is defined as follows:
$$
x^{(r)}_{k}=\overset{\infty}{\underset{i,j=1}{\sum
}}P_{ij,k}^{[0,r]}x^{(0)}_{i}x^{(0)}_{j}
$$
(iii$_{B}$)
$$
P_{ij,k}^{[s,t]}=\overset{\infty}{\underset{m,l,g,h=1}{\sum
}}P_{im,l}^{[s,r]}P_{jg,h}^{[s,r]}P_{lh,
k}^{[r,t]}x^{(s)}_{m}x^{(s)}_{g}.$$
\end{enumerate} We say that the
q.s.p. $\{P_{ij,k}^{s,t}\}$ is of {\it  type (A) or (B)} if it
satisfies the  fundamental equations (iii$_{A}$) or (iii$_{B}$),
respectively. In this definition the functions $P_{ij,k}^{[s,t]}$
denotes the probability that under the interaction of the elements
$i$ and $j$ at time $s$ the element $k$ comes into effect at time
$t$. Since for physical, chemical and biological phenomena a
certain time is necessary for the realization of an interaction,
we shall take the greatest such time to be equal to 1 (see the
Boltzmann model \cite{J} or the biological model \cite{L}). Thus
the probability $P_{ij,k}^{[s,t]}$ is defined for $t-s\geq 1$. It
should be noted that the quadratic stochastic processes are
related to quadratic transformations (see \cite{Ke},\cite{L})  in
the same way as Markov processes are related to linear
transformations.

The equations (iii$_{A}$) and  (iii$_{B}$) can be interpreted as
different laws of behavior of the "offspring".

Some examples of q.s.p. were given in \cite{G1},\cite{GM}. We note
that quadratic processes of type (A) were considered in
\cite{G1},\cite{SG1}. One of the central problems in this theory
is the study of limit behaviors of q.s.p. An ergodic principle is
a such concept relating to limit behaviors one (see \cite{K}). In
\cite{SG2} some conditions were given for q.s.p. to satisfy this
principle. It is  known \cite{G1},\cite{G2} that a certain Markov
process can be defined by means of q.s.p., therefore, it is
interesting to know the following question: if this Markov process
satisfies the ergodic principle, then would  q.s.p. satisfy that
principle? The answer to this question gives us to find more
conditions for fulfilling the ergodic principle for q.s.p., since
the theory of Markov processes are well-developed filed. In the
paper we are going to solve the formulated question for discrete
time q.s.p.

We note that a part of the results were announced in
\cite{M},\cite{SG2}.

\section{Ergodic Principle for Quadratic Stochastic Processes}

In this section we will answer to the above formulated question
for discrete time q.s.p. Before doing it, we prove some results
concerning about Markov processes.

In the sequel we will consider discrete time q.s.p., i.e. for
$\{P^{s,t}_{ij,k}\}$ the numbers $s,t$ belong to $\bn$.

Recall that a matrix $(Q_{ij})$ is called stochastic if
$$
Q_{ij}\geq 0;\ \overset{\infty}{\underset{j=1}{\sum}}Q_{ij}=1.
$$
First recall that a family of stochastic matrices
$\{(Q_{ij}^{m,n})_{i,j\in \bn}:\ m,n\in \bn, n-m\geq 1\}$ is
called  discrete time Markov process if the following condition
holds: for every $m<n<l$
\begin{equation}\label{kce}
Q_{ij}^{m,l}=\overset{\infty}{\underset{k=1}{\sum}}Q_{ik}^{m,n}Q_{kj}^{n,l}.
\end{equation}
This equation is known as the Kolmogorov-Chapman equation .

A Markov process $\{Q_{ij}^{m,n}\}$ is said to satisfy {\it the
ergodic principle} if
\begin{equation*}
\underset{n\rightarrow \infty}{\lim }|Q_{ik}^{m,n}-Q_{jk}^{m,n}|=0
\end{equation*}
is valid, for every $i,j,k,m \in \bn$. Note that this notion
firstly was introduced in \cite{K}.

 Each
$(Q_{ij}^{m,n})$-stochastic matrix defines a linear operator
${\Q}^{m,n}:\ell^{1}\rightarrow \ell^{1}$ as follows
\begin{equation}\label{mop}
({\mathbf{Q}}^{m,n}(x))_{j}=\overset{\infty}{\underset{i=1}{\sum}}Q_{ij}^{m,n}x_{i},\
\  x=(x_{n})\in \ell^{1}.
\end{equation}
Stochasticity of $(Q_{ij}^{m,n})$ implies that
\begin{equation}\label{mprop}
{\Q}^{m,n}(S)\subset S \ \ \textrm{and}\ \ \|{\Q}^{m,n}x\|_{1}\leq
\|x\|_{1}, \ \ \ x\in \ell^{1};
\end{equation}
By $\{e^{(n)}\}$ we denote standard basis of $\ell^{1}$, i.e.
$$
e^{(n)}=(\underbrace{0,0,\cdots,1}_n,\cdots),\ \ n\in \bn.
$$
Now we formulate the following well-know fact (see, for example
\cite{T}).
\begin{lem}\label{seq} A sequences $\{x_{n}\}\subset \ell^{1}$ converges
weakly if and only if it converges in norm of $\ell^{1}$.
\end{lem}We have following;
\begin{thm}\label{epm} Let $\{Q_{ij}^{m,n}\}$ be a Markov process.
The following conditions are equivalent:
\begin{enumerate}
   \item[(i)] $\{Q_{ij}^{m,n}\}$ satisfies the ergodic principle;
   \item[(ii)] For every $i,j,m \in \bn$ the  following relation holds:
$$
\lim_{n\to\infty}\|{\Q}^{m,n}e^{(i)}-{\Q}^{m,n}e^{(j)}\|_{1}=0.
$$
    \item[(iii)] For every $\varphi,\psi \in S$ and $m\in\bn$ the  following
relation holds:
$$
\lim_{n\to\infty}\|{\Q}^{m,n}\varphi-{\Q}^{m,n}\psi\|_{1}=0.
$$

\end{enumerate}
\end{thm}
\begin{pf} (i)$\Rightarrow$(ii). The ergodic principle means that the
a sequence $x_{ij,m}^{(n)}=(Q_{ik}^{m,n}-Q_{jk}^{m,n})_{k\in \bn}$
converges weakly in $\ell^{1}$, here $i,j,m\in\bn$ are fixed
numbers. According Lemma 2.1 we infer that $x_{ij,m}^{(n)}$
converges strongly in $\ell^{1}$, this means
$$
\lim_{n\to\infty}\sum_{k=1}^\infty|Q_{ik}^{m,n}-Q_{jk}^{m,n}|\rightarrow
0.
$$
On other hand, from \eqref{mop} we find
\begin{equation}\label{mo-norm}
\|{\Q}^{m,n}e^{(i)}-{\Q}^{m,n}e^{(j)}\|_{1}=\sum_{k=1}^\infty|Q_{ik}^{m,n}-Q_{jk}^{m,n}|.
\end{equation} Hence, the considered implication is
proved.

(ii)$\Rightarrow$(iii). First consider the following elements:
$\xi=\sum\limits_{i=1}^{M}\alpha_{i}e^{(i)}$,
$\eta=\sum\limits_{j=1}^{N}\beta_{j}e^{(j)}$, where
$\alpha_{i},\beta_{j}\geq 0$,
$\sum\limits_{i=1}^{M}\alpha_{i}=\sum\limits_{j=1}^{N}\beta_{j}=1$.
Using (ii) we have
\begin{eqnarray}\label{ep1}
\|{\Q}^{m,n}\xi-{\Q}^{m,n}\eta\|_{1}&= &
\|\sum\limits_{i=1}^{M}\alpha_{i}{\Q}^{m,n}e^{(i)}
-\sum\limits_{j=1}^{N}\beta_{j}{\Q}^{m,n}e^{(j)}\|\nonumber\\
&=&\|\sum\limits_{i=1}^{M}\sum\limits_{j=1}^{N}\alpha_{i}\beta_{j}{\Q}^{m,n}e^{(i)}-
\sum\limits_{i=1}^{M}\sum\limits_{i=1}^{N}\alpha_{i}\beta_{j}{\Q}^{m,n}e^{(j)}\|\nonumber\\
&\leq&
\sum\limits_{i=1}^{M}\sum\limits_{j=1}^{N}\alpha_{i}\beta_{j}
\|{\Q}^{m,n}e^{(i)}-{\Q}^{m,n}e^{(j)}\|_{1}\rightarrow 0 \ \
\textrm{as} \ \ n\to \infty.
\end{eqnarray}
Now let $\varphi,\psi \in S$ and $\varepsilon>0$. Denote
$$
G=\{\xi=\sum\limits_{i=1}^{M}\alpha_{i}e^{(i)}:\alpha_{i}\geq
0;\sum\limits_{i=1}^{M}\alpha_{i}=1,\ M\in \bn \}.
$$
It is clear that G is dense in S. Therefore, there exist $\xi$,
$\eta$ $\in G$ such that
$$
\|\varphi-\xi\|_{1}<\varepsilon/3, \ \ \
\|\psi-\eta\|_{1}<\varepsilon/3.
$$ According to \eqref{ep1} there
is $n_{0}\in \bn$ such that
$$
\|{\Q}^{m,n}\xi-{\Q}^{m,n}\eta\|_{1}<\varepsilon/3, \ \ \ \
\forall n\geq n_{0}.
$$ Hence, by means of above relations and
\eqref{mprop} we obtain
\begin{eqnarray*}
\|{\Q}^{m,n}\varphi-{\Q}^{m,n}\psi\|_{1}&\leq&
\|{\Q}^{m,n}(\varphi-\xi)\|_{1}+\|{\Q}^{m,n}(\psi-\eta)\|_{1}+\|{\Q}^{m,n}\xi-{\Q}^{m,n}\eta\|_{1}\\
&\leq
&\|\varphi-\xi\|_{1}+\|\psi-\eta\|_{1}+\varepsilon/3<\varepsilon
\end{eqnarray*}
for all $n\geq n_{0}$. Thus the implication is proved. The
implication (iii)$\Rightarrow$(i) is obvious.
\end{pf}

Let $\{P_{ij,k}^{[m,n]}\}$ be a q.s.p. Define
\begin{equation}\label{mp}
\bh_{ij}^{m,n}=\sum\limits_{l=1}^{\infty}P_{il,j}^{[m,n]}x_{l}^{(m)},
\ \ i,j \in \bn. \end{equation}
 It is clear $(\bh_{ij}^{m,n})$ is  a stochastic
matrix.
\begin{lem}\label{seq1} Let  $\{P_{ij,k}^{[m,n]}\}$ be a q.s.p. Then
$\{\bh_{ij}^{m,n}\}$ is a Markov process.
\end{lem}
\begin{pf} Consider two distinct cases with respect to types of
q.s.p.

{\tt Case (a).} Let $P_{ij,k}^{[m,n]}$ be a q.s.p of type (A).
Then we have
\begin{eqnarray*}
\sum\limits_{k=1}^{\infty}\bh_{ik}^{m,n}\bh_{kj}^{n,l}&=&
\sum\limits_{k=1}^{\infty}\left(\sum\limits_{u=1}^{\infty}P_{iu,k}^{[m,n]}x_{u}^{(m)}\right)
\left(\sum\limits_{v=1}^{\infty}P_{kv,j}^{[n,l]}x_{v}^{(n)}\right)\\
&=&
\sum\limits_{u=1}^{\infty}\left(\sum\limits_{k,v=1}^{\infty}P_{iu,k}^{[m,n]}P_{kv,j}^{[n,l]}x_{v}^{(n)}\right)x_{u}^{(m)}\\
&=&
\sum\limits_{u=1}^{\infty}P_{iu,j}^{[m,l]}x_{u}^{(m)}=\bh_{ij}^{m,l}.
\end{eqnarray*}
So $\{Q_{ij}^{m,n}\}$ is a Markov process.\\

{\tt Case (b).} Let $P_{ij,k}^{[m,n]}$ be a q.s.p  of type (B).
First consider
\begin{eqnarray*}
\sum\limits_{i,j=1}^{\infty}P_{ij,k}^{[m,n]}x_{i}^{(m)}x_{j}^{(m)}&=&
\sum\limits_{i,j=1}^{\infty}P_{ij,k}^{[m,n]}\left(\sum\limits_{e,f=1}^{\infty}P_{ef,i}^{[0,m]}x_{e}^{(0)}x_{f}^{(0)}\right)
\left(\sum\limits_{c,d=1}^{\infty}P_{cd,j}^{[0,m]}x_{c}^{(0)}x_{d}^{(0)}\right)\\
& = &
\sum\limits_{c,e=1}^{\infty}\left(\sum\limits_{f,d,i,j=1}^{\infty}P_{ef,i}^{[0,m]}
P_{cd,j}^{[0,m]}P_{ij,k}^{[m,n]}x_{f}^{(0)}x_{d}^{(0)}\right)x_{e}^{(0)}x_{c}^{(0)}\\
&=&\sum\limits_{e,c=1}^{\infty}P_{ec,k}^{[0,n]}x_{e}^{(0)}x_{c}^{(0)}=x_{k}^{(n)}.
\end{eqnarray*}
So
$$
x_{k}^{(n)}=\sum\limits_{i,j=1}^{\infty}P_{ij,k}^{[m,n]}x_{i}^{(m)}x_{j}^{(m)}.
$$
Using this equality check Markovianity:
\begin{eqnarray*}
\bh_{ij}^{m,l}&=&\sum\limits_{k=1}^{\infty}P_{ik,j}^{[m,l]}x_{k}^{(m)}=
\sum\limits_{k=1}^{\infty}\left(\sum\limits_{a,b,c,d=1}^{\infty}P_{ia,b}^{[m,n]}
P_{kc,d}^{[m,n]}P_{bd,j}^{[n,l]}x_{a}^{(m)}x_{c}^{(m)}\right)x_{k}^{(m)}\\
&=&\sum\limits_{b,d=1}^{\infty}\left(\sum\limits_{a=1}^{\infty}P_{ia,b}^{[m,n]}x_{a}^{(m)}\right)
\left(\sum\limits_{k,c=1}^{\infty}P_{kc,d}^{[m,n]}x_{k}^{(m)}x_{c}^{(m)}\right)P_{bd,j}^{[n,l]}\\
&=&\sum\limits_{b=1}^{\infty}\left(\sum\limits_{a=1}^{\infty}P_{ia,b}^{[m,n]}x_{a}^{(m)}\right)
\left(\sum\limits_{d=1}^{\infty}P_{bd,j}^{[n,l]}x_{d}^{(n)}\right)=\sum\limits_{b=1}^{\infty}\bh_{ib}^{m,n}\bh_{bj}^{n,l}.
\end{eqnarray*}
\end{pf}

The defined process $\{\bh^{m,n}_{ij}\}$ is called  {\it
associated Markov process } with respect to q.s.p. By ${\H}^{m,n}$
we denote the linear operator associated with this Markov process
(see \eqref{mop}).

We say that the {\it ergodic principle} holds for the q.s.p.
$\{P_{ij,k}^{[m,n]}\}$ if
\begin{equation*}
\underset{n\rightarrow \infty}{\lim
}|P_{ij,k}^{[m,n]}-P_{uv,k}^{[m,n]}|=0
\end{equation*}
is valid for any $i,j,u,v,k\in \mathbb{N}$ and arbitrary $m\in
\mathbb{N}$.\\
Define
\begin{equation}\label{mop2}
R_{ij}^{m,n}(x)=\sum\limits_{l=1}^{\infty}P_{il,j}^{[m,n]}x_{l},
\end{equation}
here $x=(x_{n})\in S$. It is clear that for each $m,n \in \bn$ and
$x \in S$ the matrix $(R_{ij}^{m,n}(x))$ is stochastic.
\begin{prop}\label{ep2} Let $\{P_{il,j}^{[m,n]}\}$ be a q.s.p. Then the
following conditions are equivalent:
\begin{enumerate}
   \item[(i)] For every $i,j,k \in \bn$ and $x \in S$ the following holds:
   $$
\lim_{n\to\infty}|R_{ik}^{m,n}(x)-R_{jk}^{m,n}(x)|=0.
  $$
   \item[(ii)] The Markov process $\{\bh_{ij}^{m,n}\}$ satisfies the ergodic principle.
\end{enumerate}
\end{prop}
\begin{pf}
The implication (ii)$\Rightarrow$(i). Again divide   the proof
into two cases.

{\tt Case (a).} Let $\{P_{ij,k}^{[m,n]}\}$ be a q.s.p. of type
(A). Then we have
\begin{eqnarray*}
R_{ik}^{m,n}(x)&=&\sum\limits_{l=1}^{\infty}P_{il,k}^{[m,n]}x_{l}=
\sum\limits_{l=1}^{\infty}\left(\sum\limits_{u,v=1}^{\infty}P_{il,u}^{[m,m+1]}
P_{uv,k}^{[m+1,n]}x_{v}^{(m+1)}\right)x_{l}\\
&=&\sum\limits_{l=1}^{\infty}\sum\limits_{u=1}^{\infty}
P_{il,k}^{[m,m+1]}\bh_{uk}^{m+1,n}x_{l}\\
&=&\sum\limits_{u=1}^{\infty}Q_{uk}^{m+1,n}y_{u}(i)=({\H}^{m+1,n}y(i))_{k}
\end{eqnarray*}
where
$y_{u}(i)=\sum\limits_{l=1}^{\infty}P_{il,u}^{[m,m+1]}x_{l}$.
Similarly, one gets
$$
R_{jk}^{m,n}(x)=\sum\limits_{u=1}^{\infty}\bh_{u,k}^{[m+1,n]}
y_{u}(j)=({\H}^{m+1,n}y(j))_{k}.
$$
The ergodic principle for the Markov process with Theorem
\ref{epm} implies that
$$
\|{\H}^{m+1,n}y(i)-{\H}^{m+1,n}y(j)\|_{1}\rightarrow 0\ \
\textrm{as} \ \ n\to \infty.
$$
Therefore,
\begin{eqnarray*}
|R_{ik}^{m,n}(x)-R_{jk}^{m,n}(x)|&=&|({\H}^{m+1,n}y(i))_{k}-({\H}^{m+1,n}y(j))_{k}|\\
&\leq &\|{\H}^{m+1,n}y(i)-{\H}^{m+1,n}y(j)\|_{1}\rightarrow 0\ \
\textrm{as} \ \ n\to \infty.
\end{eqnarray*}
{\tt Case (b).}  Now suppose that $(P_{ij,k}^{[m,n]})$ is a q.s.p.
of type (B). Given $x\in S$ define operator
${\mathbf{R}}^{m,n}(x):\ell^{1}\to\ell^{1}$ as follows:
$$
({\mathbf{R}}^{m,n}(x))(y)_{k}=\sum\limits_{i=1}^{\infty}R_{ik}^{m,n}(x)
y_{i},
$$
$y=(y_{i})\in \ell^{1}.$ Using stochasticity of
$(R_{ij}^{m,n}(x))$ we infer
\begin{equation}\label{norm}
\|({\R}^{m,n}(x))(y)\|_{1}\leq \|y\|_{1},\ \ \ \forall y \in
\ell^{1}.
\end{equation}
Now using (iii)$_{B}$ we find
\begin{eqnarray}\label{mp1}
R_{ik}^{m,n+1}(x)&=&\sum\limits_{l=1}^{\infty}\left(\sum\limits_{a,b,c,d=1}^{\infty}
P_{ia,b}^{[m,n]}P_{lc,d}^{[m,n]}P_{bd,k}^{[n,n+1]}x_{a}^{(m)}x_{c}^{(m)}\right)x_{l}\nonumber\\
&=&\sum\limits_{l=1}^{\infty}\sum\limits_{b,d=1}^{\infty}
\bh_{ib}^{m,n}\bh_{ld}^{m,n}P_{bd,k}^{[n,n+1]}x_{l}\nonumber\\
&=& ({\R}^{n,n+1}(y))({\H}^{m,n}(e^{(i)}))_{k},
\end{eqnarray}
here $y={\H}^{m,n}(x)$. Similarly, one gets
\begin{equation}\label{mp2}
R_{jk}^{m,n+1}(x)=({\R}^{n,n+1}(y))({\H}^{m,n}(e^{(j)}))_{k}.
\end{equation}
Therefore, it follows from \eqref{norm}-\eqref{mp2} that
\begin{eqnarray*}
|R_{ik}^{m,n+1}(x)-R_{jk}^{m,n+1}(x)|
&=&|({\R}^{n,n+1}(y))({\H}^{m,n}(e^{(i)})-{\H}^{m,n}(e^{(j)}))_{k}|\\
&\leq&
\|({\R}^{n,n+1}(y))({\H}^{m,n}(e^{(i)})-{\H}^{m,n}(e^{(j)}))\|_{1}\\
&\leq& \|{\H}^{m,n}(e^{(i)})-{\H}^{m,n}(e^{(j)})\|_{1} \rightarrow
0\ \ \textrm{as} \ \ n\to \infty.
\end{eqnarray*}
This completes the proof.
\end{pf}
\begin{prop}\label{ep3} Let $\{P_{ij,k}^{[m,n]}\}$ be a q.s.p. Then the
following conditions are equivalent:
\begin{enumerate}
   \item[(i)] $\{P_{ij,k}^{[m,n]}\}$ satisfies the ergodic principle;
   \item[(ii)] For every $x\in S$ and $i,j,k,m \in \bn $ the
   following holds:
   $$
\lim_{n\to\infty}|R_{ik}^{m,n}(x)-R_{jk}^{m,n}(x)|=0.
  $$
\end{enumerate}
\end{prop}
\begin{pf} (i)$\Rightarrow $(ii). Define a bilinear operator ${\P}^{[m,n]}:\ell^{1}\times\ell^{1}\rightarrow \ell^{1}$
as follows
$$
({\P}^{[m,n]}(x,y))_{k}=\sum\limits_{i,j=1}^{\infty}
P_{ij,k}^{[m,n]}x_{i}y_{j},
$$ where $x=(x_{n}),y=(y_{n})\in
\ell^{1}$.  According to Lemma \ref{seq} from ergodic principle we
find
$$
\lim_{n\to\infty}\|{\P}^{[m,n]}(e^{(i)},e^{(j)})-{\P}^{[m,n]}(e^{(u)},e^{(v)})\|_{1}=0
$$
for every $i,j,u,v,m \in \bn.$ The same argument of the proof of
Theorem \ref{epm} implies that
$$
\lim_{n\to\infty}\|{\P}^{[m,n]}(e^{(i)},x)-{\P}^{[m,n]}(e^{(u)},y)\|_{1}=0,
$$
for every $x,y \in S$. Hence, we have
$$
|R_{ik}^{m,n}(x)-R_{jk}^{m,n}(x)|\leq
\|{\P}^{[m,n]}(e^{(i)},x)-{\P}^{[m,n]}(e^{(j)},x)\|_{1}\rightarrow
0\ \ \textrm{as} \ \ n\to \infty.
$$
Now consider the implication (ii)$\Rightarrow$(i). From
\eqref{mop2} we have
$$
|P^{[m,n]}_{iu,k}-P^{[m,n]}_{ju,k}|=|R_{ik}^{m,n}(e^{(u)})-R_{jk}^{m,n}(e^{(u)})|\rightarrow
0\ \ \textrm{as} \ \ n\to \infty;
$$
for every $i,j,k,u \in \bn.$ Whence one gets
$$
|P^{[m,n]}_{ij,k}-P^{[m,n]}_{uv,k}|\leq
|P^{[m,n]}_{ij,k}-P^{[m,n]}_{uj,k}|+|P^{[m,n]}_{ju,k}-P^{[m,n]}_{vu,k}|
\rightarrow 0\ \ \textrm{as} \ \ n\to \infty,
$$
here we have used the equation (i) of definition q.s.p.
\end{pf}
Now we are ready to formulate our main result.
\begin{thm}\label{main} Let $\{P^{[m,n]}_{ij,k}\}$ be a q.s.p. The following
conditions are equivalent:
\begin{enumerate}
   \item[(i)] $\{P_{ij,k}^{[m,n]}\}$ satisfies the ergodic principle;
   \item[(ii)] The Markov process $\{\bh_{ikj}^{m,n}\}$ satisfies the ergodic principle.
\end{enumerate}
\end{thm}
The proof immediately follows from Propositions \ref{ep2} and
\ref{ep3}.

\section{An application of the main result}
In this section we give certain conditions for the Markov process
which ensure fulfilling the ergodic principle for q.s.p.

Now we need some auxiliary facts.
\begin{lem}\label{seq1} Let $\{a_{n}\}$ be a nonnegative sequence
which satisfies the following inequality
\begin{equation}\label{1seq}
a_{n}\leq
(1-\lambda_{n})a_{n-1}+\prod\limits_{k=1}^{n}(1-\lambda_{k}),
\end{equation}
where $\lambda_{n}\in (0,1)$,$\forall n\in\bn$ and
\begin{equation}\label{seq2}
\sum\limits_{n=1}^{\infty}\lambda_{n}=\infty,
\end{equation}
\begin{equation}\label{seq3}
\sum\limits_{j=1}^{n}\frac{\prod\limits_{k=1}^{n}(1-\lambda_{k})}{(1-\lambda_{j})}\rightarrow
0\ \ \textrm{as} \ \ n\to \infty,
\end{equation}
then $\lim\limits_{n\to\infty}a_{n}=0$.
\end{lem}
\begin{pf}
From \eqref{1seq} by iterating we get
\begin{equation}\label{seq4}
a_{n}\leq a_{1}\prod\limits_{i=2}^{n}(1-\lambda_{i})+
\sum\limits_{j=1}^{n}\frac{\prod\limits_{k=1}^{n}(1-\lambda_{k})}{(1-\lambda_{j})}.
\end{equation}
The condition \eqref{seq2} with $0<\lambda_{n}<1$ implies that
$\prod\limits_{k=1}^{n}(1-\lambda_{k})\to 0\ \ \textrm{as} \ \
n\to \infty$. Hence, by means of \eqref{seq3} we obtain $a_{n}\to
0$ as $n\to \infty$.
\end{pf}

\begin{cor}\label{seq2} Let $\{a_{n}\}$ be as above. If the sequence $\{\lambda_{n}\}$,
$(0<\lambda_{n}<1,\forall n\in\bn)$ satisfies \eqref{seq2} and the
following  relations
\begin{equation}\label{seq5}
n\prod\limits_{k=1}^{n}(1-\lambda_{k})\rightarrow 0\ \ \textrm{as}
\ \ n\to \infty
\end{equation}
\begin{equation}\label{seq6}
\sum\limits_{i=1}^{n}\frac{\lambda_{i}}{1-\lambda_{i}}=O(n)
\end{equation}
then $\lim\limits_{n\to\infty}a_{n}=0$.
\end{cor}

\begin{pf} It is enough to show that the conditions \eqref{seq5}, \eqref{seq6} imply
\eqref{seq3}. Indeed, since $0<1-\lambda_{j}<1$, we can write
$$\frac{1}{1-\lambda_{j}}=1+\varepsilon_{j},
$$ here
$\varepsilon_{j}$ is a some positive number. It then follows that
\begin{equation*}
\sum\limits_{j=1}^{n}\frac{1}{1-\lambda_{j}}=n+\sum\limits_{j=1}^{n}\varepsilon_{j}
=n+\sum\limits_{j=1}^{n}\frac{\lambda_{j}}{1-\lambda_{j}}.
\end{equation*}
It follows from  \eqref{seq6} that
\begin{equation}\label{seq7}
\sum\limits_{j=1}^{n}\frac{1}{1-\lambda_{j}}\leq C\cdot n
\end{equation} for all $n\in \bn$, here C is a constant. Therefore, from
\eqref{seq7}, \eqref{seq4} and \eqref{seq5} we infer that
\begin{eqnarray*}
\sum\limits_{j=1}^{n}\frac{\prod\limits_{k=1}^{n}(1-\lambda_{k})}{1-\lambda_{j}}&=
&\prod\limits_{k=1}^{n}(1-\lambda_{k})\sum\limits_{j=1}^{n}\frac{1}{1-\lambda_{j}}\\
&\leq& C\cdot n \prod\limits_{k=1}^{n}(1-\lambda_{k})\rightarrow
0\ \ \textrm{as} \ \ n\to \infty
\end{eqnarray*}

\end{pf}
\begin{thm}\label{epm2} Let $\{Q_{ij}^{m,n}\}$ be a Markov
process. If there exists a number $k_{0}\in \bn$ and a sequence
$\{\lambda_{n}\}$, $0<\lambda_{n}<1$,$\forall n\in\bn$ satisfying
\eqref{seq2}, \eqref{seq3} such that
\begin{equation}\label{seq8}
Q_{ik_{0}}^{n-1,n}\geq \lambda_n\ \ \textrm{for all} \ \ \ i,n\in
\bn.
\end{equation}
Then the Markov process satisfies the ergodic principle.
\end{thm}
\begin{pf} We set
$$
\sup_{i\in \bn}Q^{k,n}_{ik_{0}}=M_{k,n}(k_{0});\ \ \inf_{i\in
\bn}Q^{k,n}_{ik_{0}}=m_{k,n}(k_{0}).
$$
For $i<k<n$ we have
\begin{equation}\label{seq9}
Q^{i,n}_{ik_{0}}=\sum\limits_{l=1}^{\infty}Q^{i,k}_{il}Q^{k,n}_{lk_{0}}
\leq
M_{k,n}(k_{0})\sum\limits_{l=1}^{\infty}Q^{i,k}_{il}=M_{k,n}(k_{0}).
\end{equation}
Similarly,
\begin{equation}\label{seq10}
Q^{i,n}_{ik_{0}}\geq m_{k,n}(k_{0})
\end{equation}
By means of \eqref{seq8} we infer
\begin{equation}\label{seq11}
Q^{n-1,n}_{ik_{0}}-\lambda_{n} Q^{n-1,n}_{jk_{0}}\geq 0
\end{equation}
for all $i,j\in \bn,$ because $0\leq Q^{n-1,n}_{jk_{0}}\leq1$. It
follows
\begin{eqnarray*}
Q^{k-1,n}_{ik_{0}}&=&\sum\limits_{l=1}^{\infty}Q^{k-1,k}_{il}Q^{k,n}_{lk_{0}}\\
&=&\sum\limits_{l=1}^{\infty}[Q^{k-1,k}_{il}-\lambda_{k}Q^{k-1,k}_{jl}]Q^{k,n}_{lk_{0}}
+\lambda_{k}\sum\limits_{l=1}^{\infty}Q^{k-1,k}_{jl}Q^{k,n}_{lk_{0}}\\
&\geq&m_{k,n}(k_{0})\sum\limits_{l=1}^{\infty}[Q^{k-1,k}_{il}-\lambda
_{k}Q^{k-1,k}_{jl}]
+\lambda_{k}Q^{k-1,n}_{jk_{0}}\\
&=&(1-\lambda_{k})m_{k,n}(k_{0})+\lambda_{k}Q^{k-1,n}_{jk_{0}},
\end{eqnarray*}
whence
\begin{equation}\label{seq12}
Q^{k-1,n}_{jk_{0}}-Q^{k-1,n}_{ik_{0}}\leq
(1-\lambda_{k})(Q^{k-1,n}_{jk_{0}}-m_{k,n}(k_{0})).
\end{equation}
Since \eqref{seq12} holds for any $i,j \in \bn$, from
\eqref{seq9}-\eqref{seq10} we find
\begin{equation}\label{seq13}
M_{k-1,n}(k_{0})-m_{k-1,n}(k_{0})\leq
(1-\lambda_{k})(M_{k,n}(k_{0})-m_{k,n}(k_{0})).
\end{equation}
So iterating the last inequality we get
\begin{equation}\label{seq14}
M_{l,n}(k_{0})-m_{l,n}(k_{0})\leq
\prod\limits_{k=l+1}^{n-1}(1-\lambda_{k}).
\end{equation}
Using \eqref{kce} we have
\begin{equation*}
|Q^{m,n}_{ik}-Q^{m,n}_{jk}|\leq
\sum\limits_{l=1}^{\infty}|Q^{m,n-1}_{il}-Q^{m,n-1}_{jl}|Q^{n-1,n}_{lk},
\end{equation*}
for every  $i,j \in \bn$. Hence by means of \eqref{seq8} it yields
that
\begin{eqnarray}\label{seq15}
\sum\limits_{k=1,k\neq
k_{0}}^{\infty}|Q^{m,n}_{ik}-Q^{m,n}_{jk}|&\leq&
\sum\limits_{l,k=1}^{\infty}|Q^{m,n-1}_{il}-Q^{m,n-1}_{jl}|Q^{n-1,n}_{lk}\nonumber\\
& & -\sum\limits_{l=1}^{\infty}|Q^{m,n-1}_{il}-Q^{m,n-1}_{jl}|Q^{n-1,n}_{lk_{0}}\nonumber\\
&\leq&(1-\lambda_{n})\sum\limits_{l=1}^{\infty}|Q^{m,n-1}_{il}-Q^{m,n-1}_{jl}|.
\end{eqnarray}
Now add $|Q^{m,n}_{ik_{0}}-Q^{m,n}_{jk_{0}}|$ to both sides of
\eqref{seq15}. Then
\begin{equation*}
\sum\limits_{k=1}^{\infty}|Q^{m,n}_{ik}-Q^{m,n}_{jk}|\leq
(1-\lambda_{n})\sum\limits_{l=1}^{\infty}|Q^{m,n-1}_{il}-Q^{m,n-1}_{jl}|+
|Q^{m,n}_{ik_{0}}-Q^{m,n}_{jk_{0}}|.
\end{equation*}
Now by means of  \eqref{seq14} and \eqref{mo-norm} we infer
\begin{equation*}
\|{\Q}^{m,n}(e^{(i)})-{\Q}^{m,n}(e^{(j)})\|_{1}\leq
(1-\lambda_{n})
\|{\Q}^{m,n-1}(e^{(i)})-{\Q}^{m,n-1}(e^{(j)})\|_{1}+\prod\limits_{j=m+1}^{n-1}(1-\lambda_{j}).
\end{equation*}
Denoting $a_{n}=\|{\Q}^{m,n}(e^{(i)})-{\Q}^{m,n}(e^{(j)})\|_{1}$
and applying Lemma \ref{seq} to $a_{n}$  one gets that
\begin{equation*}
\|{\Q}^{m,n}(e^{(i)})-{\Q}^{m,n}(e^{(j)})\|_{1}\rightarrow 0\ \
\textrm{as} \ \ n\to \infty.
\end{equation*}
So, according to Theorem \ref{epm} we obtain the required
assertion.
\end{pf}

Now we can formulate the following
\begin{thm}\label{epq} Let $\{P_{ij,k}^{[m,n]}\}$ be a q.s.p.
If there exist a number $k_{0}\in \bn$ and a sequence
$\{\lambda_{n}\}$, $(0<\lambda_{n}<1)$ satisfying the conditions
\eqref{seq2}, \eqref{seq3} such that
\begin{equation}\label{seq17}
P_{ij,k_{0}}^{[n-1,n]}\geq \lambda_{n}\ \ \textrm{for all} \ \
i,j\in \bn,
\end{equation}
then  the q.s.p. satisfies the ergodic principle.
\end{thm}
\begin{pf} Consider the Markov process $\{\bh_{ij}^{m,n}\}$
associated with given q.s.p. (see \eqref{mp}). Then \eqref{seq17}
implies that
\begin{equation*}
\bh_{ik_{0}}^{n-1,n}=\sum\limits_{l=1}^{\infty}P_{il,k_{0}}^{[n-1,n]}x^{(n-1)}_{l}
\geq
\lambda_{n}\sum\limits_{l=1}^{\infty}x^{(n-1)}_{l}=\lambda_{n}.
\end{equation*}
Consequently, the Markov process satisfies the conditions of
Theorem \ref{epm2}. Therefore the ergodic principle is valid for
it. Now by means of Theorem \ref{main} we infer that q.s.p.
satisfies the ergodic principle.
\end{pf}

\section*{acknowledgements}

The third named author (F.M.) thanks NATO-TUBITAK for providing
financial support and Harran University for kind hospitality and
providing all facilities. The work is also partially supported by
Grant $\Phi$-1.1.2 of Rep. Uzb.

\end{document}